\documentclass[
final
]{dmtcs-episciences}

\newcommand{\set}[1]{\left\{ #1 \right\}}
\newcommand{\Ord}{{\cal O}}
\newcommand{\kl}[1]{\left(\, #1 \,\right)}
\newcommand{\const}{\sf const}
\newcommand{\nn}{\mathbb{N}_{_{0}}}
\newcommand{\nat}{\mathbb{N}}
\newcommand{\pinprim}{p \in \mathbb{P}}
\newcommand{\prim}{\mathbb{P}}
\newcommand*\teilt[2]{#1 \! \left| #2 \right.}
\newcommand{\ejk}{1 \leq j \leq k}
\newcommand{\B}{{\Box}}
\newcommand{\njk}{0 \leq j \leq k}
\newcommand*\Betr[1]{\left|#1\right|}

\newtheorem{kor}{Corollary}[section]
\newtheorem{lemma}{Lemma}[section]
\newtheorem{satz}{Theorem}[section]
\newtheorem{remark}{Remark}[section]
\newtheorem{bsp}{Example}[section]

\usepackage{amsmath}

\usepackage[utf8]{inputenc}
\usepackage{subfigure}

\usepackage[round]{natbib}

\author{Alexander Büchel\affiliationmark{1,2}
  \and Ulrich Gilleßen\affiliationmark{2}
  \and Kurt-Ulrich Witt\affiliationmark{1,2}}
\title[An output-sensitive Algorithm to partition a Sequence of Integers into Subsets with equal Sums]{An output-sensitive Algorithm to partition a\\ 
 Sequence of Integers into Subsets with equal Sums}
\affiliation{
  b-it Applied Science Institute, University Bonn-Rhein-Sieg, Sankt Augustin, Germany\\
  Research Group Discrete Mathematics and Optimization (ADIMO), University Bonn-Rhein-Sieg, Sankt Augustin, Germany}
\received{2017-07-18}
\revised{2018-07-18, 2018-11-27}
\accepted{2018-12-04}
\begin{document}
\publicationdetails{20}{2018}{2}{18}{4994}
\maketitle
\begin{abstract}
We present a polynomial time algorithm, which solves a nonstand\-ard variation of the well-known PARTITION problem: Given positive integers $n, k$ and $t$ such that $t \geq n$ and $k \cdot t = {n+1 \choose 2}$, the algorithm partitions the elements of the set $I_n = \set{1, \ldots, n}$ into $k$ mutually disjoint subsets $T_j$ such that $\cup_{j=1}^k T_j = I_n$ and $\sum_{x \in T_{j}} x = t$ for each $j \in \set{1,2, \ldots, k}$. The algorithm needs $\Ord${\tiny $\kl{\! n \! \cdot \!  \kl{\! \frac{n}{2k} \! + \! \log \frac{n(n+1)}{2k} \! } \! }$} steps to insert the $n$ elements of $I_n$ into the $k$ sets $T_j$
\end{abstract}

\section{Introduction}

For 
\begin{math}
	n \in \mathbb{N}
\end{math} 
let
\begin{math}
	I_n = \set{1, \ldots, n}
\end{math}
be the set of integers from $1$ to $n$, and 
\begin{math}
	\Delta_n = \frac{n(n+1)}{2}
\end{math}
the sum of these elements. In this paper we consider a variant of the PARTITION problem and present a solution for a class of special instances of this variant. The general version of our variant is given by 
\begin{math}
	n, k, t_1, \ldots, t_k \in \mathbb{N},
\end{math}
and the question is whether there exists $k$ pairwise disjoint subsets 
\begin{math}
	T_j \subseteq I_n
\end{math}
such that the elements of $T_j$ add up to $t_j$, and the union of these sets equals $I_n$. We call such a collection of sets $T_j$ a 
\begin{math}
	(t_1,t_2, \ldots, t_k)
\end{math}
\!-partition of $I_n$. 

\medskip

\cite{Fu:1992} show, that for
\begin{math}
	k, l, t \in \mathbb{N}
\end{math}
with 
\begin{math}
	0 < l \leq \Delta_{n}
\end{math}
and 
\begin{math}
	(k-1)t + l + \Delta_{k-2} = \Delta_{n}
\end{math}
a 
\begin{math}
	(t, t+1, \ldots, t+k-2, l)
\end{math}
\!-partition of $I_n$ exists. \cite{Chen:2005} prove, that a 
\begin{math}
	(t_{1}, \ldots, t_{k})
\end{math}
\!-partition of $I_n$ exists, if
\begin{math}
	\sum_{j=1}^{k} t_{j} = \Delta_{n}
\end{math}
and $t_{j} \geq t_{j+1}$ for $1 \leq j \leq k-1$ and $t_{k-1} \geq n$ hold. In \cite{Buechel:2016} we present a $0$/$1$-linear program to solve partition problems.

\medskip

In the special case, where $t_j = t = \const$ we call $T_1, \ldots, T_k$ a $(k,t)$-partition of $I_n$. Given $n, k, t \in \mathbb{N}$ with $t \geq n$ and $\Delta_n = k \cdot t$ the decision problem reduces to the question, whether a $(k,t)$ partition of $I_n$ exists. \cite{Straight:1979} show that for all $k, t$ with $\Delta_{n} = k \cdot t$ and $t \geq n$ a partition of $I_n$ exists.  \cite{Ando:1990} withdraw the condition $\Delta_{n} = k \cdot t$ and prove that for positive integers $n$, $k$ and $t$, the set $I_n$ contains $k$ disjoint subsets having a constant sum $t$ if and only if 
\begin{math}
	k(2k-1) \leq k \cdot t \leq \Delta_n.
\end{math}

\medskip

Where as the cited papers study for which $k$-tuples $(t_1, \ldots, t_k)$-partitions of $I_n$ exist, we are interested in efficient algorithms to determine partitions. In this paper we consider problem instances $\Pi(n,k,t)$ with $t \geq n$ and $\Delta_n = k \cdot t$. In Section \ref{sectionAlgorithm} we introduce the recursive algorithm $\Pi\mathit{Solve}$ which determines a partition for each instance $\Pi(n,k,t)$. Before, in Section \ref{sectionMeander} we present the so called meander algorithm which solves problem instances $\Pi(n,k,t)$, where $n$ is even and $2k$ is a divisor of $n$ or where $n$ is odd and $2k$ divides $n+1$, respectively. The reason is, that $\Pi\mathit{Solve}$ can be stopped, when one of these conditions is reached, and the remaining partition can be determined directly by means of the meander algorithm. In Section \ref{sectionComplexity} we analyze the run time complexity of $\Pi\mathit{Solve}$. Section \ref{sectionConclusion} summarizes the paper and mentions some ideas to improve $\Pi\mathit{Solve}$. 

\medskip

Inputs for the algorithms are $n$, $k$ and $t$, hence these have length $\Ord(\log n)$. Since it is to be expected that the complexity to insert $n$ elements into $k$ sets is at least $\Ord(n)$, we will consider the complexity of the algorithms not depending on the size of the inputs, but output-sensitive, i.e. depending on $n$ and $k$.

\section{Meander Algorithm}\label{sectionMeander}

For $a \in \nn$ and $b \in \nat$ we denote $\teilt{b}{a}$ if $b$ is a divisor of $a$. Given the problem instance $\Pi(n,k,t)$ the meander algorithm applies if $n$ is even and $\teilt{2k}{n}$ or if $n$ is odd an $\teilt{2k}{n+1}$, respectively. The algorithm distributes the elements of the set $I_n$ into the subsets $T_j$ such that these sets build a $(k,t)$-partition of $I_n$, i.e. the sets $T_j$ fulfill the conditions
\begin{align}
T_i \cap T_j & = \emptyset, \,\, 1 \leq i, j \leq k, \ i \not= j \label{glMeander102} \\ 
\bigcup_{j=1}^k T_j & = I_n \label{glMeander101} \\ 
\sum_{x \in T_j} x & = t, \,\, \ejk \label{glMeander1}
\end{align}

\subsection{Case: n even and 2k$\mid$\!n}\label{subsectionnger}

%\subsection{Case: \begin{math}\mathbf{n}\end{math} even and \begin{math}\mathbf{\teilt{2k}{n}}\end{math}} \label{subsectionnger}

Figure \ref{figMeanderEven} shows the part of the meander algorithm which solves problem instances $\Pi(n,k,t)$ when $n$ is even and $2k$ divides $n$. To prove that the algorithm determines a correct $(k,t)$-par\-ti\-tion of $I_n$ we have to show that the 
%set of elements assigned to the subsets $T_j$ in (1) and (2) is equal to $I_n$ and (ii) that the resulting partition fulfills condition (\ref{glMeander1}). 
partition fulfills the conditions above. Condition (\ref{glMeander102}) is obviously fulfilled.  We will verify (\ref{glMeander101}) in Lemma \ref{lemmaMeandereven1} and (\ref{glMeander1}) in Lemma \ref{lemmaMeandereven2}. 

\medskip
Let
\begin{align}
X_1(n,k) & = \set{2ki - (j-1) \mid 1 \leq i \leq \frac{n}{2k}, \, \ejk } \\
X_2(n,k) & = \set{2k(i-1) + j \mid 1 \leq i \leq \frac{n}{2k}, \, \ejk } 
\end{align}
be the sets of elements of $I_n$ which are distributed in assignment (I) or assignment (II), respectively.

\medskip

\begin{figure}[h]
	\rule{14.8cm}{1pt} 
	
	\verb| |\texttt{\textbf{meandereven}$(n, k, t)$;} \\ 
	\verb| |\texttt{\textbf{input:} $\phantom{i}n, k, t$ with $n$ even, $\teilt{2k}{n}$, $t \geq n$, and $\Delta_n = k \cdot t$;} \\
	\verb| |\texttt{\textbf{output:} $(k,t)$-partition $T_j, \ejk$, of $I_n$;} \\
	\verb|    |\texttt{$T_j := \emptyset, \,\, 1 \leq j \leq k$;} \\
	\verb|    |\texttt{\textbf{for}} $j := 1$ \textbf{to} $k$ \textbf{do} \\
	\verb|       |\texttt{\textbf{for}} $i := 1$ \textbf{to} $\frac{n}{2k}$ \textbf{do} \\
	\verb|          |\texttt{ (I) $T_j : = T_j \cup \set{2ki-(j-1)}$}; \\ 
	\verb|          |\texttt{(II) $T_j := T_j \cup \set{2k(i-1) + j)}$}; \\ 
	\verb|       |\texttt{\textbf{endfor}}; \\
	\verb|    |\texttt{\textbf{endfor}}; \\
	\verb| |\texttt{\textbf{end.}}
	
	\rule{14.8cm}{.5pt} 
	\textbf{\caption{\label{figMeanderEven}Meander Algorithm in case $\boldmath{n}$ even and $\boldmath{\teilt{2k}{n}}$.}}
	
	\rule{14.8cm}{1pt} 
\end{figure}

{\lemma\label{lemmaMeandereven1} \, Let $\Pi(n,k,t)$ be a problem instance such that $n$ even and $\teilt{2k}{n}$, then 
\begin{math}
	I_n = X_1(n,k) \cup X_2(n,k).
\end{math}}
 
\begin{proof}
For each $x \in I_n$ there exist unambiguously $i, r$ such that 
\begin{align}
x & = 2k(i-1) + r, \, 1 \leq i \leq \frac{n}{2k}, \, 1 \leq r \leq 2k \label{glMeandereven3}
\end{align}
We consider the two following sets of remainders $r \in I_{2k}$: 
\begin{math}
	R_1 = \set{2k - (j-1) \mid \ejk}
\end{math}
and 
\begin{math}
	R_2 = \set{j \mid \ejk}.
\end{math}
Since $r \in R_1$, if $k+1 \leq r \leq 2k$, it follows $R_1 \cap R_2 = \emptyset$ and $R_1 \cup R_2 = I_{2k}$. Thus with respect to (\ref{glMeandereven3}) we get either 
\begin{equation}
x = 2k(i-1) + 2k - (j-1) = 2ki - (j-1)
\end{equation}
or
\begin{equation}
x = 2k(i-1) + j 
\end{equation}
It follows 
\begin{math}
	x \in X_1(n,k) \cup X_2(n,k).
\end{math}
Hence we have shown 
\begin{math}
	I_n \subseteq  X_1(n,k) \cup X_2(n,k).
\end{math}

\medskip

If $x \in X_1(n,k)$, then $k+1 \leq x \leq n$, and if $x \in X_2(n,k)$ then $1 \leq x \leq n-k$. Thus, if 
\begin{math}
	x \in X_1(n,k) \cup X_2(n,k),
\end{math} 
we have $1 \leq x \leq n$, hence $x \in I_n$ and thereby 
\begin{math}
	X_1(n,k) \cup X_2(n,k) \subseteq I_n.
\end{math}
\end{proof}

{\lemma\label{lemmaMeandereven2} \, Let $\Pi(n,k,t)$ be a problem instance with $n$ even and $\teilt{2k}{n}$, then the output $T_j$, $\ejk$, of $\mathtt{meandereven}(n,k,t)$ fulfills condition (\ref{glMeander1}).}

\begin{proof}
For each $j \in \set{1, \ldots, k}$ we have: $\sum_{x \in T_j} x = $
\begin{align*}
\sum_{i=1}^{\frac{n}{2k}} (2ki - (j-1)) + \sum_{i=1}^{\frac{n}{2k}} (2k(i-1) + j)
 = 2k \sum_{i=1}^{\frac{n}{2k}} (2i-1) + \frac{n}{2k} = 2k \frac{n^2}{4k^2} + \frac{n}{2k}
 = \frac{n(n+1)}{2k} = t
\end{align*} 
\end{proof}

{\satz\label{satMeandereven} \, $\mathtt{meandereven}(n,k,t)$ 

\medskip

\textbf{a)} determines a correct partition of $I_n$ for all problem instances $\Pi(n,k,t)$ with $n$ even and $\teilt{2k}{n}$,  

\medskip

and

\medskip

\textbf{b)} needs $\Ord(n)$ steps to insert the $n$ elements of $I_n$ into the sets $T_j$.}

\begin{proof} 
a) follows immediately from Lemmas \ref{lemmaMeandereven1} and \ref{lemmaMeandereven2}, and b) is obvious.
\end{proof}

\subsection{Case: n odd and 2k$\mid$n+1}\label{subsectionnuger}

%\subsection{Case: $\boldmath{n}$ odd and $\boldmath{\teilt{2k}{n \! + \!  1}}$}\label{subsectionnuger}

To solve problem instances $\Pi(n,k,t)$ with $n$ odd and $\teilt{2k}{n+1}$ we adapt slightly the $\mathtt{meandereven}$-algorithm (see Fig. \ref{figMeanderOdd}). The correctness of the $\mathtt{meanderodd}$-algorithm can be shown analogously to the proof of the correctness of the $\mathtt{meandereven}$-algorithm. At this point we define the sets of elements assigned due to labels (I) and (II) in the $\mathtt{meanderodd}$-algorithm as
\begin{align} 
X_1'(n,k) & = \set{2ki - j \mid  1 \leq i \leq \frac{n+1}{2k}, \, \ejk } \\ 
X_2'(n,k) & = \set{2k(i-1) + (j-1) \mid 1 \leq i \leq \frac{n+1}{2k}, \, \ejk }
\end{align}

\medskip

\begin{figure}[h]
	\rule{14.8cm}{1pt} 
	
	\verb| |\texttt{\textbf{meanderodd}$(n, k, t)$;} \\ 
	\verb| |\texttt{\textbf{input:} $\phantom{i}n, k, t$ with $n$ odd, $\teilt{2k}{n+1}$, $t \geq n$, and $\Delta_n = k \cdot t$;} \\
	\verb| |\texttt{\textbf{output:} $(k,t)$-partition $T_j, \ejk$, of $I_n$;} \\
		\verb|    |\texttt{$T_j := \emptyset, \,\, \ejk$;} \\
			\verb|    |\texttt{\textbf{for}} $j := 1$ \textbf{to} $k$ \textbf{do} \\
	\verb|       |\texttt{\textbf{for}} $i := 1$ \textbf{to} $\frac{n}{2k}$ \textbf{do} \\
	\verb|          |\texttt{ (I) $T_j := T_j \cup \set{2ki - j}$}; \\ 
	\verb|          |\texttt{(II) $T_j := T_j \cup \set{2k(i-1) + (j-1)}$}; \\ 
	\verb|       |\texttt{\textbf{endfor}}; \\
	\verb|    |\texttt{\textbf{endfor}}; \\
	\verb| |\texttt{\textbf{end.}}
	
	\rule{14.8cm}{.5pt} 
	\textbf{\caption{\label{figMeanderOdd}Meander Algorithm in case $\boldmath{n}$ odd and $\boldmath{\teilt{2k}{n+1}}$.}}
	
	\rule{14.8cm}{1pt} 
\end{figure}

\medskip

{\remark \, In order to avoid a case distinction, we first assign the element $0$ ($i=1$, $j=1$) to set $T_1$. For this reason, in the following we assume that $I_n$ contains the element $0$, too.} \hfill $\B$

{\lemma\label{lemmaMeanderodd1} \, Let $\Pi(n,k,t)$ be a problem instance such that $n$ odd and $\teilt{2k}{n+1}$, then 
\begin{math}
	I_n = X_1'(n,k) \cup X_2'(n,k).
\end{math}} 
	
\begin{proof}
For each $x \in I_n$ there exist unambiguously $i, r$ such that 
\begin{align}
	x & = 2k(i-1) + r, \, 1 \leq i \leq \frac{n+1}{2k}, \, 0 \leq r \leq 2k - 1 \label{glMeanderodd3}
\end{align}
We consider the sets of remainders $r \in I_{2k-1}$: 
\begin{math}
	R_1' = \set{2k - j \mid \ejk}
\end{math}
and 
\begin{math}
	R_2' = \set{j-1 \mid \ejk} 
\end{math}
\begin{math}
	= \set{j \mid \njk-1}.
\end{math}
Since $r \in R_1'$, if $k \leq r \leq 2k-1$, it follows $R_1 \cap R_2 = \emptyset$ and $R_1 \cup R_2 = I_{2k-1}$. Thus with respect to (\ref{glMeanderodd3}) we get  
\begin{equation}
	x = 2k(i-1) + 2k - j = 2ki - j 
\end{equation}
or
\begin{equation}
	x = 2k(i-1) + (j-1)
\end{equation}
respectively. It follows 
\begin{math}
	x \in X_1'(n,k) \cup X_2'(n,k).
\end{math}
Hence we have shown
\begin{math}
	I_n \subseteq  X_1'(n,k) \cup X_2'(n,k).
\end{math}
	
\medskip
	
If $x \in X_1'(n,k)$, then $k \leq x \leq n$, and if $x \in X_2'(n,k)$ then $0 \leq x \leq n-k$. Thus, if 
\begin{math}
	x \in X_1'(n,k) \cup X_2'(n,k),
\end{math}
we have $0 \leq x \leq n$, hence $x \in I_n$ and thereby 
\begin{math}
	X_1'(n,k) \cup X_2'(n,k) \subseteq I_n.
\end{math}
\end{proof}

{\lemma\label{lemmaMeanderodd2} \, Let $\Pi(n,k,t)$ be a problem instance with $n$ odd and $\teilt{2k}{n+1}$, then the output $T_j$, $\ejk$, of $\mathtt{meanderodd}(n,k,t)$ fullfills condition (\ref{glMeander1}). }

\begin{proof}	
For each $j \in \set{1, \ldots, k}$ we have 
\begin{align*}
	\sum_{x \in T_j} x & = 
	\sum_{i=1}^{\frac{n+1}{2k}} (2ki - j) + \sum_{i=1}^{\frac{n+1}{2k}} (2k(i-1) + (j-1)) \\ 
	& = 2k  \sum_{i=1}^{\frac{n}{2k}} (2i-1) - \frac{n+1}{2k} = 2k \frac{(n+1)^2}{4k^2} - \frac{n+1}{2k} \\
	& = \frac{n(n+1)}{2k} = t
\end{align*} 
\end{proof}

{\satz\label{satMeanderodd} \, $\mathtt{meanderodd}(n,k,t)$ 
	
\medskip

\textbf{ a)} determines a correct partition of $I_n$ for all problem instances $\Pi(n,k,t)$ with $n$ odd and $\teilt{2k}{n+1}$, 

\medskip

and

\medskip

\textbf{ b)} needs $\Ord(n)$ steps to insert the $n$ elements of $I_n$ into the sets $T_j$.}

\begin{proof}
a) follows from Lemmas \ref{lemmaMeanderodd1} and \ref{lemmaMeanderodd2}, and b) is obvious.
\end{proof}
		
%\section{The Algorithm $\boldmath{\Pi\mathit{Solve}}$}\label{sectionAlgorithm}
\section{The Algorithm $\Pi$\emph{Solve}}\label{sectionAlgorithm}

In this section we present the different cases which the $\Pi\mathit{Solve}$-algorithm distinguishes using ideas similar to those used in \cite{Straight:1979}. The input to the algorithm are the integers $n, k, t \in \nat$ with $t \geq n$ and $\Delta_n = k \cdot t$. The output is a $(k,t)$-partition $T_j$, $\ejk$, of $I_n$, which fullfills condition (\ref{glMeander1}). We prove that the algorithm works correctly in all cases.
	
%\subsection{Case: $\boldmath{2n > t}$}\label{subsection2n1lt}
\subsection{Case: 2n $>$ t}\label{subsection2n1lt}
		
In this case the algorithm makes a distinction between the cases $t$ even and $t$ odd. 
		
\subsubsection{Case: t even} 
		
The algorithm starts with filling $\frac{2n-t}{2}$ sets as follows:
\begin{align}
	T_j & = \set{ t - n + (j-1), n - (j-1) }, \, 1 \leq j \leq \frac{2n-t}{2} \label{glteven1}
\end{align}
Obviously these sets are disjoint and fullfill condition (\ref{glMeander1}). The union of these sets is the set \linebreak 
\begin{math}
	\set{t-n, \ldots, \frac{t}{2} - 1, \frac{t}{2} + 1, \ldots, n}.
\end{math}
Thus the elements of the set $I_{t-n-1}$ %$1, \ldots, t-n-1$ 
and the element $\frac{t}{2}$ remain, these have to be distributed into the empty $k - \frac{2n-t}{2}$ sets. To do this, each of these sets is split into two subsets:
\begin{align}
	T_j & = T_{j,1} \cup T_{j,2}, \, \frac{2n-t}{2} + 1 \leq j \leq k \label{glteven2}
\end{align}
The total number of these subsets is $2(k-n) + t$. The set $T_{\frac{2n-t}{2}+1,1}$ is filled with the element $\frac{t}{2}$:
\begin{align}
	T_{\frac{2n-t}{2}+1,1} & = \set{ \frac{t}{2} } \label{glteven3}
\end{align}
Thus it remains to distribute the elements of $I_{t-n-1}$ into the $2(k-n) + t - 1$ sets $T_{\frac{2n-t}{2}+1,2}$ and $T_{j,s}$, $\frac{2n-t}{2} + 2 \leq j \leq k$, $s \in \set{1,2}$, i.e. it remains to solve the problem instance $\Pi(n',k',t')$ where
\begin{align}
	n' & = t - n - 1 \label{glns1}\\ 
	k' & = 2(k-n) + t - 1 \label{glks1}\\
	t' & = \frac{t}{2} \label{glts1}
\end{align}
We have to verify that this instance fulfills the input conditions 
\begin{align}
	\Delta_{n'} & = k' \cdot t' \label{glveri03}
\end{align}
and
\begin{align}
	t' & \geq n' \label{glveri031}
\end{align}
Using (\ref{glns1}) -- (\ref{glts1}) we get on one side
\begin{align}
	\Delta_{n'} & = \frac{n'(n'+1)}{2}  = \frac{(t-n-1)(t-n)}{2} = \Delta_n + \frac{t^2 - 2tn - t}{2} \label{glveri13}
\end{align}
and on the other side
\begin{align}
	k' \cdot t' & = (2(k - n) + t - 1) \cdot \frac{t}{2} = k \cdot t + \frac{t^2 -2tn - t}{2} \label{glveri23} 
\end{align}
Since for our initial problem $\Pi(n,k,t)$ the condition $\Delta_n = k \cdot t$
holds, the verification of (\ref{glveri03}) follows immediately from (\ref{glveri13}) and (\ref{glveri23}). 

\medskip

From $2n > t$ immediately follows $\frac{t}{2} > t - n - 1$. Using (\ref{glns1}) and (\ref{glts1}) condition (\ref{glveri031}) is verified, too.

\medskip

Thus the algorithm can recursively continue to solve the initial problem by determining a solution for the instance $\Pi(n',k',t')$. 

\subsubsection{Case: t odd} 
		
In this case the algorithm initially fills $\frac{2n-t+1}{2}$ sets as follows:
\begin{align}
	T_j & = \set{ t - n + (j-1), n - (j-1) }, \, 1 \leq j \leq \frac{2n-t+1}{2} \label{gltodd1}
\end{align}
Obviously these sets are disjoint and fullfill condition (\ref{glMeander1}). 
The union of these sets builds the set \linebreak $\set{t-n, \ldots, n}$. Thus the elements of the set $I_{t-n-1}$ %$1, \ldots, t-n-1$ 
remain, these have to be distributed into the empty $k - \frac{2n-t+1}{2}$ sets. Therefore, the instance $\Pi(n',k',t')$ has to be solved, where
\begin{align}
	n' & = t - n - 1  \label{glns2} \\ 
	k' & = k - \frac{2n-t+1}{2} \label{glks2} \\
	t' & = t \label{glts2} 
\end{align}
To proof that this instance is feasible we have to verify, that the input conditions (\ref{glveri03}) and (\ref{glveri031}) are fulfilled in this case as well.

\medskip

Using (\ref{glns2}) -- (\ref{glts2}) we get on one side	 
\begin{align}
	\Delta_{n'} & = \frac{n'(n'+1)}{2} = \frac{(t-n-1)(t-n)}{2} = \Delta_n + \frac{t^2 - 2tn - t}{2} \label{glveri11}
\end{align}
and on the other side
\begin{align}
	k' \cdot t' & = \kl{ k - \frac{2n - t + 1}{2} } \cdot t = k \cdot t + \frac{t^2 -  2tn - t}{2} \label{glveri21} 
\end{align}
Since $\Delta_n = k \cdot t$ the verification of (\ref{glveri03}) follows immediately from (\ref{glveri11}) and (\ref{glveri21}).

\medskip

From $2n > t$ it follows $n > t - n - 1$. From this we get by means of the input condition $t \geq n$ and the definitions (\ref{glns2}) and (\ref{glts2}): 
\begin{math}
	t' = t \geq n > t - n - 1 = n',
\end{math}
i.e. condition (\ref{glveri031}) is fulfilled. 
		
\subsection{Case: 2n $\leq$ t}\label{subsection2nlesst}
		
In this case each set $T_j$ is split into two disjoint subsets: 
\begin{math}
	T_j = T_{j,1} \cup T_{j,2}, \, \ejk.
\end{math}
The sets $T_{j,1}$ will be filled as follows:
\begin{align}
	T_{j,1} & = \set{ n - 2k + j, n - (j-1) } \label{gl2nt1}
\end{align}
Hence the elements $n-2k+1, \ldots, n$ are already distributed, and the two elements in  each of these sets add up to
\begin{align}
	n - (i-1) + n - 2k + i & = 2(n-k) + 1
\end{align}
It remains to partition the elements of $I_{n-2k}$ 
%$1, \ldots, n-2k$ 
into the sets $T_{j,2}$ such that the sum of elements in each $T_{j,2}$ equals $t - (2(n-k)+1)$. Thus it remains to solve the problem instance $\Pi(n',k',t')$ with
\begin{align}
	n' & = n - 2k \label{gl219} \\ 
	k' & = k \label{gl220} \\
	t' & = t - 2(n-k) - 1 \label{glt2k} 
\end{align}
As well as in the former cases we have to assure, that the input conditions (\ref{glveri03}) and (\ref{glveri031}) are fulfilled.
%\begin{align}
%	\Delta_{n'} & = k' \cdot t' \label{glveri02}
%\end{align}
On the one side we have
\begin{align}
	\Delta_{n'} & = \frac{(n-2k)(n-2k+1)}{2} = \Delta_n + 2k^2 - k - 2kn \label{glveri12}
\end{align}
and on the other side
\begin{align}
	k' \cdot t' & = k \cdot \kl{ t - 2(n-k) - 1 } = k \cdot t -2kn + 2k^2 - k \label{glveri22}
\end{align}
(\ref{glveri03}) follows immediately from (\ref{glveri12}) and (\ref{glveri22}).%, and $\Pi\mathit{Solve}$ can be recursively applied to $\Pi(n',k',t')$.

\medskip

From $t \geq 2n$ it follows $n + 1 \geq 4k$. By subtraction we get 
\begin{math}
	t - n - 1 \geq 2n - 4k
\end{math}
and from this and definitions (\ref{gl219}) and (\ref{glt2k})  
\begin{math}
	t' = t - 2n + 2k - 1 \geq n - 2k = n', 
\end{math} 
i.e. condition (\ref{glveri031}) is verified.

\medskip

The considerations so far lead to the algorithm $\Pi\mathit{Solve}$ shown in Figure \ref{figure1}, and we proved that it works correctly in all cases.

\begin{figure}[p]
	\rule{14.8cm}{1pt}

	\verb| |\texttt{$\boldmath{\Pi\mathit{Solve}(n,k,t)}$;} \\
	\verb| |\texttt{\textbf{ input:}  $\phantom{i}n, k, t$ with $t \geq n$, and $\Delta_n = k \cdot t$;} \\
	\verb| |\texttt{\textbf{ output:} $(k,t)$-partition $T_j, \ejk$, of $I_n$;} \\
	
	\medskip
	
	\hspace*{0.5cm}  (I) \texttt{\textbf{ case} $\teilt{2k}{n}$} \\ \\ 
	\hspace*{1.05cm} \texttt{\textbf{ then} fill $\set{T_j}_{1 \leq j \leq k}$ by $\mathit{meandereven}(n,k,t)$} \\ \\ 
	\hspace*{0.8cm} \texttt{\textbf{ case} $\teilt{2k}{n+1}$} \\ \\ 
	\hspace*{1.05cm} \texttt{\textbf{ then} fill $\set{T_j}_{1 \leq j \leq k}$ by $\mathit{meanderodd}(n,k,t)$} \\ \\ 
	\hspace*{0.5cm} (II) \texttt{\textbf{ case} $t \geq 2n$} \\ \\ 
	\hspace*{1.05cm} \texttt{\textbf{ then} \textbf{ for} $1 \leq j \leq k$ \textbf{ do} $T_{j,1} = \set{n-2k+j, n-(j-1)}$ \textbf{ endfor};} \\ \\ 
	\hspace*{2.25cm} \texttt{fill $\set{T_{j,2}}_{1 \leq j \leq k}$ by $\Pi\mathit{Solve}(n-2k,k,t-2(n-k)-1))$;} \\ \\ 
	\hspace*{2.25cm} \texttt{\textbf{ for} $1 \leq j \leq k$ \textbf{ do} $T_{j} = T_{j,1} \cup T_{j,2}$ \textbf{ endfor}} \\ \\ 
	\hspace*{0.5cm}(III) \texttt{\textbf{ case} $t < 2n$ \textbf{ and} $t$ even \\ \\ 
	\hspace*{1.05cm} \textbf{ then} \textbf{ for} $1 \leq j \leq \frac{2n-t}{2}$ \textbf{ do} $T_{j} = \set{t-n+(j-1), n-(j-1)}$ \textbf{ endfor};} \\ \\ 
	\hspace*{2.25cm} \texttt{$T_{\frac{2n-t}{2}+1,1} = \set{\frac{t}{2}}$;} \\ \\ 
	\hspace*{2.25cm} \texttt{fill $T_{\frac{2n-t}{2}+1,2}$, $\set{T_{j,1}}_{\frac{2n-t}{2}+2 \leq j \leq k}$ and $\set{T_{j,2}}_{\frac{2n-t}{2}+2 \leq j \leq k}$} \\ \\ 
	\hspace*{2.55cm} \texttt{by $\Pi\mathit{Solve}(t-n-1,2(k-n)+t-1,\frac{t}{2})$;} \\ \\
	\hspace*{2.25cm} \texttt{\textbf{ for} $\frac{2n-t}{2} + 1 \leq j \leq k$ \textbf{ do} $T_{j} = T_{\frac{2n-t}{2} + j,1} \cup T_{\frac{2n-t}{2} + j,2}$ \textbf{ endfor}} \\ \\  
	\hspace*{0.5cm}(IV) \texttt{\textbf{ case} $t < 2n$ \textbf{ and} $t$ odd \\ \\ 
	\hspace*{0.9cm} \textbf{ then for} $1 \leq j \leq \frac{2n-t+1}{2}$ \textbf{ do} $T_{j} = \set{t-n+(j-1), n-(j-1)}$ \textbf{ endfor};} \\ \\ 
	\hspace*{2.25cm} \texttt{fill $\set{T_{j}}_{\frac{2n-t+1}{2}+1 \leq j \leq k}$ by $\Pi\mathit{Solve}(t-n-1,k-\frac{2n-t+1}{2},t)$} \\ \\ 
	\texttt{\textbf{ end.}}
	
	\rule{14.8cm}{.5pt} 
	\textbf{ \caption{\label{figure1}Algorithm $\boldmath{\Pi\mathit{Solve}}$.}}
	
	\rule{14.8cm}{1pt} 
\end{figure}

\medskip

\section{Complexity}\label{sectionComplexity}

In this section we analyse the worst case run time complexity of the $\Pi\mathit{Solve}$-Algorithm. The algorithm consists of four subalgorithms related to the cases we distinguish: (I) $\teilt{2k}{n}$ or $\teilt{2k}{n+1}$,  (II) $t \geq 2n$, (III) $t < 2n$ and $t$ even, (IV) $t < 2n$ and $t$ odd. We abbreviate these cases by $m$ (meander), $s$ (smaller), $ge$ (greater even), and $go$ (greater odd),  respectively. Then the run $\Pi\mathit{Solve}(n,k,t)$ can be represented by a sequence 
\begin{math}
	\rho'(n,k,t) \in \set{m, s, ge, go}^{+}\!.
\end{math}

{\bsp \, \textbf{ a)} Let $n = 1337$. The list of runs for all partitions of $I_{1337}$ is:
\begin{align*}
	\rho'(1337,3,298151) & = m \\ 
	\rho'(1337,7,127779) & = s^{94}\,\mathit{go}\,m \\ 
	\rho'(1337,21,42593) & = s^{30}\,\mathit{go}\,s\,\mathit{ge}\,m \\ 
	\rho'(1337,191,4683) & = ss\,\mathit{go}\,m \\ 
	\rho'(1337,223,4011) & = m \\ 
	\rho'(1337,573,1561) & = \mathit{go}\,m \\ 
	\rho'(1337,669,1337) & = m 
\end{align*}
\textbf{ b)} Let $n = 9999$, then we have 
\begin{align*}
	\rho'(9999,4444,11250) & = \mathit{ge}\,s^3\,\mathit{ge}^4\,\mathit{go}\,m \\ 
	\rho'(9999,4040,12375) & = \mathit{go}\,s^4\,\mathit{go}\,s^4\,\mathit{go}\,s\,\mathit{ge}\,m \\ 
	\rho'(9999,3960,12625) & = \mathit{go}\,s^{3}\,\mathit{ge}\,\mathit{go}\,s^{8}\,\mathit{go}\,m \\ 
	\rho'(9999,3333,15000) & = \mathit{ge}^3\,\mathit{go}\,m \\ 
	\rho'(9999,12,4166250) & = s^{415}\,\mathit{go}\,s\,\mathit{ge}^2\,m
\end{align*}} %\hfill $\B$

\medskip

Let $\alpha$ be a non empty sequence over $\Omega' = \set{m, s, ge, go}$, then $\mathit{first}(\alpha)$ is the first and $\mathit{last}(\alpha)$ the last symbol of $\alpha \in {\Omega'}^+$, and $\mathit{head}(\alpha)$ is the sequence without the last symbol. $\Betr{w}_a$ is the number of occurrences of symbol $a \in \Omega'$ in the sequence $w \in \Omega'^\ast$.

\medskip

%\begin{minipage}{13cm}

Obviously we have

{\lemma\label{lemmam} \, Let $\Pi(n,k,t)$ be a problem instance, then $\mathit{last}(\rho'(n,k,t)) = m$ and $m$ is not a member of $\mathit{head}(\rho'(n,k,t))$. } \hfill$\B$

%\end{minipage}

\medskip

Thus, we may neglect the last symbol of $\rho'(n,k,t)$ and denote $\rho(n,k,t) = \mathit{head}(\rho'(n,k,t))$. As well we do not need the alphabet $\Omega'$, because $\rho(n,k,t) \in \set{s, \mathit{ge}, \mathit{go}}^{\ast}$. We denote this alphabet by $\Omega$.

\medskip

Next we show, that the last call before the recursion stops with the $m$-case cannot be $s$.

{\lemma\label{lemmasnotlast}  \, Let $\Pi(n,k,t)$ be a problem instance. If $\Betr{\rho(n,k,t)} \geq 1$, then $\mathit{last}(\rho(n,k,t)) \not= s$.}

\begin{proof}
We assume $\mathit{last}(\rho(n,k,t)) = s$. Let $\Pi(\nu, \kappa, \tau)$ be the problem instance before the last $s$-call. Then by (\ref{gl219}) and (\ref{gl220}) after the $s$-call we have $\nu' = \nu - 2\kappa$ and $\kappa' = \kappa$. Since the next call is $m$ it has to be $\teilt{2\kappa'}{\nu'}$ or $\teilt{2\kappa'}{\nu' + 1}$, thus we have $\teilt{2\kappa}{\nu - 2\kappa}$ or $\teilt{2\kappa}{\nu - 2k + 1}$. It follows $\teilt{2\kappa}{\nu}$ or $\teilt{2\kappa}{\nu + 1}$. Hence the instance $\Pi(\nu,\kappa, \tau)$ would have been solved by an $m$-call, a contradiction to our assumption $\mathit{last}(\rho(n,k,t)) = s$.
\end{proof}

{\kor\label{kornotlast} \,%Let $\Pi(n,k,t)$ be a problem instance such that 
If $\Betr{\rho(n,k,t)} \geq 1$, then 
\begin{math}
	\mathit{last}(\rho(n,k,t)) \in \set{\mathit{ge}, \mathit{go}}.
\end{math} } \hfill $\B$

\subsection{Case: 2n $>$ t and t odd}

From $2n > t$ we can conclude $t > 2(t - n - 1)$. Using (\ref{glns2}) and (\ref{glts2}) we get $t' > 2n'$. This leads to

{\lemma\label{lemmagos} \, Let $\Pi(n,k,t)$ be a problem instance with $2n > t$, $t$ odd and 
\begin{math}
	\rho'(n,k,t) = \alpha \mathit{go} \beta, \alpha \in \Omega^{\ast}\!, \linebreak \beta \in \Omega'^+\!,
\end{math}
then 
	
\medskip

\textbf{ a)} $\mathit{first}(\beta) = m$, if $\Betr{\beta} = 1$,
	
\medskip

\textbf{ b)} $\mathit{first}(\beta) = s$, if $\Betr{\beta} \geq 2$. } \hfill$\B$

\medskip

Thus, after the case $\mathit{go}$ the recursion ends by call of the meander algorithm or the recursion continues with the $\mathit{s}$ case either.

{\kor\label{korgos} \, Let $\Pi(n,k,t)$ be a problem instance with $2n > t$ and $t$ odd, then 
\begin{equation}
\Betr{\rho(n,k,t)}_s \geq \Betr{\rho(n,k,t)}_{go}.
\end{equation} } %\hfill $\B$

\subsection{Case: 2n $>$ t and t even}

From (\ref{glts1}) it follows immediately

\begin{align}
	\Betr{\rho(n,k,t)}_{ge} & \leq \log t = \log \frac{n(n+1)}{2k} \label{gllogt}
\end{align}

\subsection{Case: 2n $\leq$ t}\label{subsection2nlesstc}

In this case if the algorithm performs the instance $\Pi(n,k,t)$, then the next instance to solve may be $\Pi(n',k,t')$ with $n' = n - 2k$ and $t' = t - 2(n-k)-1$ (cf. Subsection \ref{subsection2nlesst}, equations (\ref{gl219}) and (\ref{glt2k}), respectively). By $n^{(\ell)}$ and $t^{(\ell)}$ we denote the value of $n$ and $t$ in the $\ell^{\text{th}}$ recursion call in the case $2 n^{(\ell)} \leq t^{(\ell)}$. Thus we have 
\begin{math}
	n^{(0)} = n, \, n^{(1)} = n' = n - 2k
\end{math}
and 
\begin{math}
	t^{(0)} = t, \, t^{(1)} = t' = t - 2(n-k) - 1,
\end{math}
for example. By induction we get
\begin{align}
	n^{(\ell)} & = n - 2k \cdot \ell \label{gl2kell} \\ 
	t^{(\ell)} & = t - 2n \cdot \ell + 2k \cdot \ell^2 - \ell \nonumber \\
			   & = t - (2(n - k \cdot \ell) + 1) \cdot \ell \label{gl2tell}
\end{align}
Now we determine the order of the maximum value of $\ell$ guaranteeing the condition $2 n^{(\ell)} \leq t^{(\ell)}$. Using (\ref{gl2kell}) and (\ref{gl2tell}) we get
\begin{align}
	0 & \leq t^{(\ell)} - 2 n^{(\ell)} \\
	  & = t - (2(n - k \cdot \ell) + 1) \cdot \ell - 2(n - 2k \cdot \ell)
\end{align}
To determine $\ell$ we solve the quadratic equation
\begin{align}
	0 & = \ell^2 + \frac{4k - 2n - 1}{2k} \cdot \ell + \frac{t - 2n}{2k}
\end{align}
which has the solutions
\begin{align}
	\ell_{1,2} & = - \frac{4k - 2n - 1}{4k} \pm \sqrt{\kl{\frac{4k - 2n - 1}{4k}}^2 - \frac{t - 2n}{2k}} \\ \\ 
			   & = - \frac{4k - 2n - 1}{4k} \pm \frac{4k - 1}{4k}
\end{align}
i.e.
\begin{align}
	\ell_1 & = \frac{n}{2k}, \,\,\, \ell_2 = \frac{n+1}{2k} - 2
\end{align}
Finally we get
\begin{align}
	\ell & \leq \frac{n}{2k}
\end{align}
Thus, we have just proven

{\lemma\label{lemma2kkt} \,Let $\Pi(n,k,t)$ be a problem instance. If $\rho(n,k,t) = s^\ell x$ with $x \in \set{\mathit{ge}, \mathit{go}}$, then $\ell \leq \frac{n}{2k}$. }\hfill$\B$

\medskip

\begin{minipage}{13cm}

Corollary \ref{korgos}, inequality (\ref{gllogt}) and Lemma \ref{lemma2kkt} lead to

{\satz\label{satzconclusion} \, Let $\Pi(n,k,t)$ be a problem instance.
	
\medskip

\textbf{ a)} \, Then the recursion depth of $\Pi\mathit{Solve}(n,k,t)$ is $\Ord\kl{\frac{n}{2k} + \log \frac{n(n+1)}{2k}}$.

\medskip

\textbf{ b)} \, Since the complexity of operations the algorithm performs in each recursion call (assig\-ning elements of $I_n$ to some set $T_j$, arithmetic comparisons and operations) is $\Ord(n)$ it follows that $\Pi\mathit{Solve}$ needs 
\begin{equation}
	\Ord\kl{n \cdot \kl{\frac{n}{2k} + \log \frac{n(n+1)}{2k}}}
\end{equation}
steps to insert the $n$ elements of $I_n$ into the $k$ sets $T_j$. \hfill$\B$ }

\end{minipage}

\section{Conclusion}\label{sectionConclusion}
		
In Section \ref{sectionAlgorithm} we present the recursive algorithm $\Pi\mathit{Solve}$ which solves following special PARTITION problems $\Pi(n,k,t)$: Given $n, k, t \in \nat$ with $t \geq n$ and $\Delta_n = k \cdot t$, then the algorithm partitions the set $I_n = \set{1, \ldots, n}$ into $k$ mutually disjoint sets such that the elements in each set add up to $t$. The recursion can be stopped, if $n$ ist even and ${2k}$ is a divisor ${n}$ or if $n$ is odd and ${2k}$ is a divisor of ${n+1}$, respectively, because in these cases the meander algorithms presented in Section \ref{sectionMeander} can be applied, which directly determines a partition.

\medskip

We prove that the algorithm works correctly and needs 
\begin{equation}
	\Ord\kl{n \cdot \kl{\frac{n}{2k} + \log \frac{n(n+1)}{2k}}}
\end{equation}
steps to assign the elements of $I_n$ to the $k$ subsets $T_j$ for each problem instance $\Pi(n,k,t)$. Taking into account that the algorithm for the inputs $n$ and $k$ determines an output consisting of $k$ sets to which the elements of $I_n$ are to be distributed so that all constraints are met, $\Pi\mathit{Solve}$ is a polynomial output-sensitive time algorithm.

\medskip

In \cite{Jagadish:2015} an approximation algorithm for the cutting sticks-problem is presented. Because the cutting sticks-problem can be transformed into an equivalent partitioning problem our algorithms can be applied to the corresponding cutting sticks-problems.

\medskip

Further research may investigate whether ideas from the previous chapters and cited papers can be used to improve the efficiency of the $\Pi\mathit{Solve}$-algorithm. In \cite{Buechel:2016}, \cite{Buechel:2017a} and \cite{Buechel:2017b} we present efficient solutions for problem instances $\Pi(n,k,t)$, where $n = q \cdot k$, $q, k$ odd; $n = m^2 - 1$, $m \geq 3$; $n = p - 1$, $n = p$, $n = 2p$, $\pinprim$, where $\prim$ is the set of prime numbers. Thus we may augment the $\Pi\mathit{Solve}$-algorithm by related conditions to stop further recursion calls. 

\medskip

\acknowledgements

We would like to thank Arkadiusz Zarychta, a member of the ADIMO group as well, who created a tool by means of which we are able to test the algorithm and to analyse experimentally its performance. 

\medskip

Furthermore we would like to thank the reviewers for their valuable comments leading to improvements of the presentations.

\nocite{*}
\bibliographystyle{abbrvnat}
% use the following instead if you encounter problems 
%\bibliographystyle{alpha}
\bibliography{lit}
\label{sec:biblio}

\end{document}